\documentclass[prl,twocolumn,showpacs,preprintnumbers,amsmath,amssymb]{revtex4}
\usepackage{latexsym}
\usepackage{graphicx}
\usepackage{pstricks,pst-node,pst-tree}
\usepackage{dcolumn}
\usepackage{bm}
\usepackage{xcolor}
\usepackage{float}
\usepackage{epsfig}
\input amssym.def
\input amssym.tex

\newcommand{\bea}{\begin{eqnarray}}
\newcommand{\eea}{\end{eqnarray}}

\newcommand{\nn}{\nonumber}

\begin{document}
\title{ Heavy-Tail Distribution from Correlation of Discrete Stochastic Process }


\author{Jongwook Kim${}{}^{\natural}\footnote{Electronic correspondence:~\bf{jwkim@apctp.org}}$ and Teppei Okumura${}{}^{\dagger}$
}

\affiliation{{~}\\
\mbox{${}^{\natural}$Asia Pacific Center for Theoretical Physics, Pohang 790-784, Korea.
}
\\
\mbox{${}^{\dagger}$Institute for the Early Universe, Ewha Womans University, Seoul 120-750, Korea.
}
}

\begin{abstract}
We propose a stochastic process driven by the memory effect with novel distributions which include both exponential and leptokurtic heavy-tailed distributions. A class of the distributions is analytically derived from the continuum limit of the discrete binary process with the {\it renormalized auto-correlation}\,. The moment generating function with a closed form is obtained, thus the cumulants are calculated and shown to be convergent. The other class of the distributions is numerically investigated. The combination of the two stochastic processes of memory with different signs under {\it regime switching mechanism} does result in behaviors of power-law decay. Therefore we claim that {\it memory is the alternative origin of heavy-tail}.  
\end{abstract}

\maketitle

{\it Introduction}\,. - Non-Gaussian distributions can be observed everywhere in nature. They are prevalent in condensed matter systems and soft matter systems such as crystal growth, polymer transportation, and polymer distribution \cite{Sornette,deGennes}. They are also observed in astrophysics, for instance, galaxy mass and velocity distributions \cite{PressSchechter, Peebles}. Leptokurtosis, skewness, and the power-law tail are  characteristics of the social science data such as financial time series  \cite{Kleinert,Cont} and complexities \cite{Barabasi}\,.
Among the non-Gaussianities, the most notable features are the power-law distribution and the auto-correlation, which are often referred to as heavy-tail and memory. 

{\it Heavy-tail} is the statistical phenomenon which predicts extreme values more frequently than the conventional Gaussian law. The Levy model \cite{Kleinert, Cont, Rachev}\,, derived from the Poisson jump diffusion process, has been used to describe this feature. The Levy distribution demonstrates excellent fits to various kinds of real data. Moreover the jump diffusion is Markovian, hence the Levy model has an advantage in the use of Ito calculus \cite{Ito}. However, this Markovian non-Gaussian model cannot explain non-Markovian behaviors. {\it Memory} is the representative non-Markovian quantity which is often measured by the auto-correlation of a stochastic process. Burst and clustering are the typical phenomena which are driven by memory. Recently the heavy-tail statistics  in terms of correlation has been studied in the context of physics \cite{Krug}. Fractional Brownian motion \cite{Mandelbrot} (FBM) is one of the most popular mathematical theories in dealing with the memory effect. The parameter of FBM is called the Hurst coefficient \cite{Hurst} which is defined as the exponent of the power-law correlation function of the two white noises at different times. While FBM has an obvious advantage on the issue of memory, there are three drawbacks. First, conventional stochastic calculus is not applicable to FBM due to the very non-Markovian nature, therefore little is known for its analytic properties. Second the inner dynamics of the Hurst parameter is unclear. Third, heavy-tail distribution cannot be constructed by FBM. 

In this letter, we propose a new micro process driven by memory that results in the heavy-tail distribution. This non-Markovian non-Gaussian process not only goes beyond the Poisson jump diffusion process, which is the only known heavy-tail micro process so far, but also overcomes the three drawbacks of FBM in that it is an analytically-solvable model.

\begin{figure}
\begin{center}
\includegraphics[width=.48\textwidth]{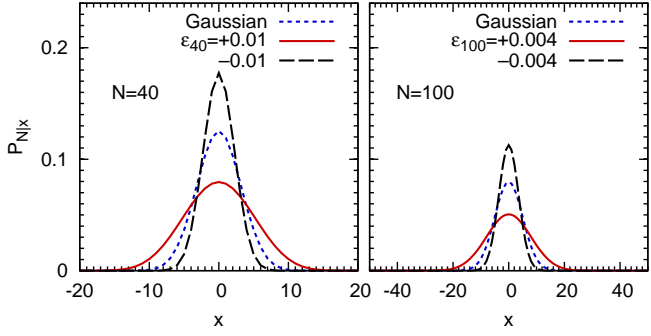}
\caption{Probability distribution functions $P_{n|x}$ defined in equation (\ref{recurrence}) with 
$N=40$ (left) and $N=100$ (right). The values of $\epsilon$ are chosen to have $\epsilon=\pm 0.4/N$.
The Gaussian distribution where $\epsilon=0$ is also shown for comparison.
}
\label{fig1}
\end{center}
\end{figure}

{\it Stochastic process with memory}\,. - We design the simplest possible stochastic process with both  non-Markovian and non-Gaussian properties. For such purpose, we implement a binomial process and  assume zero-skewness for simplicity. We require that the single recurrence rule should define the whole process and that the number of parameters should be minimal. 
Under these restrictions there are not many possibilities to model the stochastic process. We propose a process \footnote{It is a game of tossing a flexible coin that remembers the previous outcome and bends itself by the result at each process. For example, when the outcome is turned out to be the head it bents itself toward the tail in a certain degree and the process will be positively auto-correlated. The size of these bents are denoted as $\epsilon$ which is going to be matched to the size of auto-correlation in the stochastic process.}  that is defined by the recursive formula for the discrete probability distribution function (PDF) $P_{n|x}$ with binomial {\it transfer probabilities} given as
 \footnote{This stochastic process does not go in the class of  Urn-process implied by the $\beta$-distribution \cite{FriedmanUrn} though it has close resemblance. Novel feature is the use of the renormalized auto-correlation, which enables one to arrive at the convergent theory.},
\bea
P_{n|x}&=&P_{n-1|x+1} \cdot [{1\over 2} - {(x+1) }{\epsilon}]  \nn 
\\
	& &~~~ +	P_{n-1|x-1} \cdot [{1\over 2} + {(x-1 )}{\epsilon }]\, \label{recurrence}
\eea
where $n\in\{0,1,2,\cdots,N\}$ and we define stochastic displacements from the origin by the accumulation of the  $+1$s and the $-1$s. Each event is indexed as $n$ for the entire $N$ events. The displacement $x$ runs from $-n$ to $n$ in steps of $2$ {\it i.e.} $x \in \{ -n\,, -n + 2\,, -n+4 \,, \cdots\,, n-2\,, n\}$ with $n+1$ elements and we start with $P_{0|0} = 1$. For $x$ which is not in this range $P_{n|x}=0$. Probability distributions for some values of the coupling $\epsilon$ are illustrated in Fig. \ref{fig1}\,. We introduce the generating function for the analytic study of the process as
$Z_n(q)=\sum_{{x\over 2}=-{n\over 2}}^{n\over 2} P_{n|x} q^x 
				\equiv {_{2}\sum}_{x=-n}^{n} P_{n|x} q^x \,,$ where $Z_n(1)= {_{2}\sum}_{x=-n}^{n} P_{n|x} = 1$\,.
The distribution is obtained simply by reading off $q^2$ coefficients from the generating function\,, and the boundary conditions are $Z_0(q)=1$ and $Z_1(q)={1\over 2}(q+q^{-1})\, \label{BCs}$\,.
In the Gaussian limit the moment generating function is given by $
\lim_{\epsilon \rightarrow 0} Z_N(q) = {1\over 2^N} (q+q^{-1})^N \,.
$
Then the recurrence (\ref{recurrence}) is recast into a differential equation for $q$, 
\bea
Z_{n+1}(q) = {1\over 2}(q+q^{-1}) Z_n(q)
		+ {\epsilon }(q^2-1) \partial_q Z_n(q) \,. \label{Zequationlinear_x}
\eea
We make a substitution that $Z_n(q) = \epsilon^n q^{1\over 2\epsilon}(q^2 - 1)^{-{1\over 2\epsilon}}Y_n(q) $\, for all $n$\,. Then 
\bea
Y_{n+1}(q) = 
			(q^2-1) \partial_q Y_n(q) \,,
\eea
with $Y_0(q) = (q -{1\over q})^{1\over 2\epsilon}$\,.
Further we substitute
\bea
r={1\over 2}\ln{1-q \over 1+q} = -\tanh^{-1}(q)\,,
\eea
then $q={1-e^{2r} \over 1+ e^{2r}}=-\tanh(r)$\,, $Y_0(q(r)) 
 				= \Big({\sinh(2r) \over 2}\Big)^{-{1\over 2\epsilon}}$\,. In this coordinate, we have
\bea
Z_N(q(r)) = \epsilon^N \Big({\sinh(2r)\over 2}\Big)^{1\over 2\epsilon} [\partial_r]^N \Big({\sinh(2r)\over 2}\Big)^{-{1\over 2\epsilon}} \,.
\eea
Solving $N$-th derivative, the generating function is obtained in the closed form as\,,
\bea
	& & Z_N(q) = (2\epsilon)^{N-1} 
 	 {{{1\over 2\epsilon}+N}\choose{N}}\times~\nn \\
 	 & & 
	\sum_{k=1}^N 
	\frac
	{
		{1\over 4}^{k}  {{N}\choose{k}} \sum_{i=0}^k (-1)^{i+k} {{k}\choose{i}} (2i -k)^N ({1-q\over 1+q})^{2i}
		(q+{1\over q}+2)^k
	}
	{{1\over 2\epsilon}+k}\,, \label{Zq}\nn \\
\eea
where the binomial function is defined as the product
${{{1\over 2\epsilon}+N}\choose{N} } = {\prod_{i=0}^{N-1}({1\over 2\epsilon}+N-i)/ N!}\,.$

\begin{figure}
\begin{center}
\includegraphics[width=.48\textwidth]{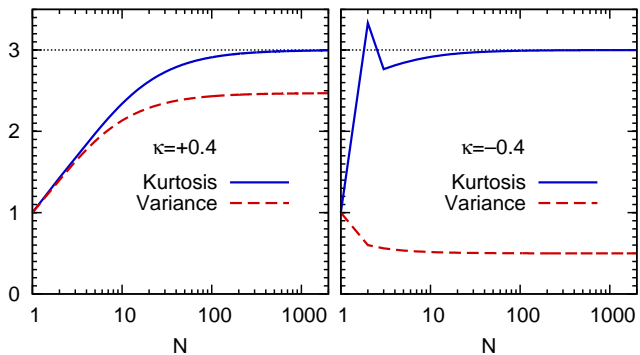}
\caption{Variance and kurtosis as a function of $N$ with $\kappa=0.4$ (left) and $\kappa=-0.4$ (right).
The values of the variance are divided by $N$. 
}
\label{fig2}
\end{center}
\end{figure}

{\it Convergence : Cumulants and Auto-Correlation}\,. - The coupling $\epsilon$ is a fixed constant at all times. We introduce the renormalization of the coupling $\epsilon \rightarrow \epsilon_N={\kappa\over N}$, hence the more fine-grained the time lattice is, the smaller the coupling becomes. By this renormalization, the convergence is guaranteed (as one can infer from Fig. \ref{fig1}) and the interpretation of the auto-correlation formula (\ref{autocorrelation}) given below becomes natural. Instead of $\epsilon$, we take $\kappa$ as the control parameter of the model, which ranges from $-{1\over 2}$ to ${1\over 2}$. 
Cumulants or $n$-point correlations are calculated from derivatives of the generating function (\ref{Zq})\,. 
As expected, the average and the skewness are zero. The variance and the $4$-point correlation are calculated to be
\bea
\langle x^2 \rangle 
&=&(2\epsilon)^{N-1}{{{1\over 2\epsilon}+N}\choose{N}}\times \label{var}\\  
& &				\sum_{k=1}^N {(-1)^k\, k \over{1\over 2\epsilon}+k}{{N}\choose{k}}
					\, {1 \over 2}[-(2-k)^N + (-k)^N]\,, \nn \\
^{\epsilon\rightarrow 0}&\rightarrow& N + 4{N(N-1)\over 2!}\cdot \epsilon \nn \\
& &~~~+ {16}{N(N-1)(N-2)\over 3!}\cdot \epsilon^2 + \mathcal{O}(\epsilon^3)\,. \nn
\eea
\bea
\langle x^4 \rangle 
&=& (2\epsilon)^{N-1}{{{1\over 2\epsilon}+N}\choose{N}} 
			\sum_{k=1}^N (-1)^{k+1}k \epsilon{{N}\choose{k}} \times\label{kur} \\ 
& & ~~~\Big\{-4(2-k)^N + 3(4-k)^N +(-k)^N +\nn \\
& &3k[2(2-k)^N - (4-k)^N -(-k)^N]\Big\}/{(2+4k\epsilon)}\nn \\
^{\epsilon\rightarrow 0}&\rightarrow& N(3N-2) + 8{N(N-1)\over 2!}(3N-4)\cdot \epsilon \nn \\
& & +56\, {N(N-1)(N-2) \over 3!}(3N -{43\over7})\cdot \epsilon^2 + \mathcal{O}(\epsilon^3)\,.  \nn 
\eea
Now we present the analytic evidence of the convergence by showing the convergence of cumulants at the large $N$ limit. It is calculated analytically in the small $\kappa$ region and for the convergence at high $\kappa$, we resort to numerical simulation. 
Taking the coupling $\epsilon$ as the renormalized one ${\kappa \over N}$ and making use of the results (\ref{var}) and (\ref{kur}), we calculate the kurtosis up to second order in $\kappa$. Nontrivial cancellations occur in each coefficient of $\kappa$ and $\kappa^2$. Consequently in the large $N$ limit the kurtosis converges to $3$ for the sufficiently small $\kappa$. 
\bea
\lim_{N \to \infty}{\langle x^4 \rangle\over (\langle x^2 \rangle)^2} = 3 + 0\cdot \kappa + 0\cdot \kappa^2 +\mathcal{O}(\kappa^3) \,. \label{kur_kappa}
\eea
To confirm the convergence, we numerically calculate the kurtosis. The measured kurtosis is shown in Fig. \ref{fig2} as a function of $N$, corroborating the analytic result (\ref{kur_kappa}) at the large $N$ limit. 
Here we present the size of auto-correlation with generic time lag $L-1$ with $L= 2\,,3\,,\cdots$. The displacement at the level $n+1$ is denoted as $\Delta S_{n+1}$ to have value $\pm 1$ and $\Delta S_{n+L}$ to have alternating values of $\pm 1$, but they are weighted by the probabilities defined by (\ref{recurrence}) . The result is
\bea
& & E[\Delta S_{n+1} \Delta S_{n+L}] \nn \\
&=& 2 \epsilon_N (1+ 2\epsilon_N)^{L-2} (1 + 2\epsilon_N \langle x^2 \rangle) \label{autocorrelation}\\
&=& {2 \kappa \over N} (1+ 2 H \kappa) + \mathcal{O}({1\over N^2})\,, \nn
\eea
where $ H\equiv (1+2\kappa + {8\over 3}\kappa^2 + \mathcal{O}(\kappa^3) )$\,. Note that $n$ and $x$ parametrize time and displacement. The auto-correlation has no $n$ dependence therefore the process is stationary for changes in time. Moreover a stochastic trajectory at large displacement is shown to be highly correlated, which is consistent with our intuition. 
The second equality not only exhibits the convergence for the auto-correlation but also clarifies the origin of the Hurst coefficients \cite{Hurst} whose inner dynamics had been unclear. The auto-correlation is well defined in the large $N$ limit, when the result (\ref{var}) is incorporated. The variance numerically calculated in Fig. \ref{fig2} is also consistent with Eq. (\ref{var})\,. As briefly mentioned in the introduction, the Hurst coefficient $h$ is defined from the displacement in terms of the exponent of the time lag {\it i.e.} $\Delta x \equiv 2 \langle \sqrt{x^2}\rangle  \sim 2 N^{h}$. In our theory $h$ is read from $\Delta x \sim 2 \sqrt{N \cdot{H}}$, and the value of $h$ is simply ${1/2}$ in the Gaussian case since $ \Delta x \sim 2 \sqrt N$\,. In Eq. (\ref{autocorrelation}), the auto-correlation formula is recast in terms of $H$, which means that the Hurst coefficient $h$ can be determined from the micro dynamics of the stochastic process. It is one of the novel features in our theory that differs from the conventional formulae based on FBM \cite{Mandelbrot}.

{\it Heavy-tail : Regime switching }\footnote{J.K. is thankful to Gabjin Oh of his suggestion for the stochastic process with thresholds.}\,. - We construct the leptokurtic-heavy-tail distribution profile from the memory effect. The PDF of a process with negative $\kappa$ has a sharp peak in the center but falls to zero very rapidly, while the PDF with positive $\kappa$ has slowly falling tails but a lower peak in the center than that with negative $\kappa$, as illustrated in Fig. \ref{fig1}. 

Regardless of signs of the auto-correlation, the kurtosis converges to 3 in case of the monotonic memory process defined by (\ref{recurrence}). Such inflexibility is not desirable when we consider general applications. Inspired by {\it regime switching mechanism}\,\cite{RSM}\,, we combine two stochastic processes with different signs of the autocorrelation\,, that allows us to arrive at the heavy-tail distribution with flexible values of kurtosis. As an example, we examine a regime switching behavior of the Gaussian profile with $x$ dependence such as, $\epsilon \rightarrow \epsilon ( b-e^{-{x^2 / 2 \delta^2}} )\,$ with $b>0$. Under this  mechanism, a specific random walk flips the sign of its auto-correlation to be positive in the occasion of deviation from the center regime thereby the trail tends to accelerate to the marginal displacement. Remaining in the center regime on the other hand, the trajectory is negatively correlated and tends to concentrates increasingly to the center. 
\begin{figure}
\includegraphics[width=.48\textwidth]{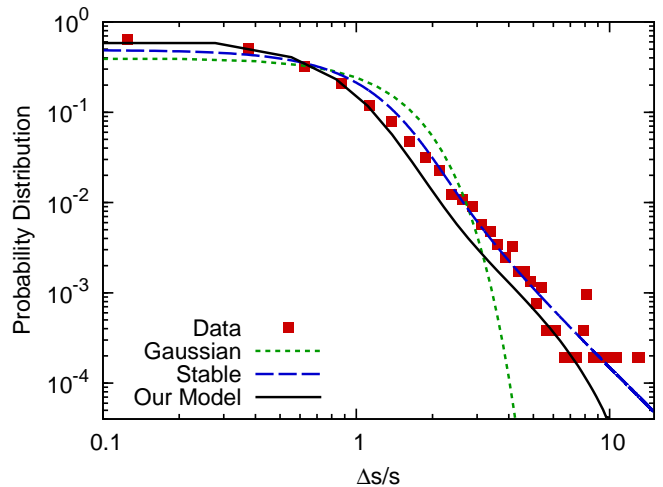}
\caption{The probability distribution of log return data of the Dow Jones daily index. 
The dotted and dashed lines are predictions from Gaussian and stable distributions respectively. 
The solid line is our model with the best fit parameters, $b=0.38$, $\delta=13.8$ and $\kappa=3.0$. 
Note that the PDFs of the data and predictions from Gaussian distribution and our model are scaled 
to $\sigma=1$, while that from stable distribution is scaled so that fitting to the data becomes best. 
}
\label{fig3}
\end{figure}
Making use of this model, we show a proper quantitative match to the Dow Jones data \footnote{ Log return of the Dow Jones daily index from 1928.10.1
to 2012.1.23 shown as dotted line in the figure}. In Fig. \ref{fig3} we illustrate the log-scaled probability distribution of the log return, 
$\Delta s/s = (s_{i+1} - s_{i})/s_i$, where $s_i \equiv s(t_i)$  is the stock price at a given day $t_i$. 
The PDF is normalized to have the variance $\sigma=1$. 
Obviously the data cannot be explained by the Gaussian distribution which is depicted as the dotted line. In comparing the stock data with our model prediction, we calculate $\chi^2$ statistics for logarithmic values of the PDF by adopting $b$, $\delta$ and $\kappa$ as free parameters. 
Since we are interested in the heavy-tail behaviors which starts at $\Delta s/s \approx 2$, 
the fit is made for the range $1.5\leq \Delta s/s \leq 10$. Our results are, however, not largely 
affected by this choice. The black solid line shown in Fig. \ref{fig3} is our prediction with 
the best fit parameters, $(b,\delta,\kappa)=(0.38,13.8,3.0)$. Unlike the monotonic memory process (\ref{recurrence}), $\kappa$ can be large as long as the transfer probability is between 0 and 1 \footnote{In other words, the process is well defined as long as the inequality $(x-1)\kappa / N( b-Exp({-{x^2 / 2 \delta^2}}) ) < 1/2$\, holds.}. In this data analysis, the location of $x=\Delta s/s$ where the transfer probability starts to be ill-defined is around $\Delta s/s=7$. The majority of the data points fall in the well defined region.

We also compare the performance of our model with the existent heavy-tail theory {\it i.e.} Levy theory. The blue dashed line is the $\alpha$-stable distribution. At a glance the stable distribution seems to have a better fit to the heavy-tail data, but there are several problems. Stable distributions generically have divergent cumulants due to its strict power-law behavior at high $\Delta s/s$ and it is impossible to compute the variance, which is one of the most important statistics. Therefore we have no choice but to treat the scaling parameter as a free parameter, which means we choose the scaling such that we necessarily obtain the best fit. Nevertheless, stable distribution fails to explain the PDF data at low $\Delta s/s$. On the contrary, we normalize the variance to be unity for our model to compare fairly to the data, and it precisely predicts the data at low $\Delta s/s$. A tempered-stable distribution \cite{Rachev} can be used to cure the defect of a stable distribution, by artificially assuming an exponentially-decaying regulator distribution function at the asymptotic region. On the other hand our model can naturally comprise such a regulating behavior from the physical mechanism {\it i.e.} regime switching. 

{\it Discussion}\,. - 
From a simple micro process with the memory effect we have derived a class of analytically solvable distributions of good flexibility controlled by the single auto-correlation parameter $\kappa$. The resulting distribution has an exponential tail but it can be made {\it heavier} around the intermediate region. The PDF with a closed form has to be calculated by performing an inverse Laplace transform or fast Fourier transform to the moment generating function (\ref{Zq}), and this will be presented in future works. We studied an application to the real financial data so as to demonstrate the accuracy of our model. Issues on convergence and renormalization should be clarified for the regime switching model with further analytic study, however our numerical analysis provides an important first step.
 
Both non-Markovian (FBM) and Markovian non-Gaussian (Levy) models have been extensively developed, however the non-Markovian non-Gaussian framework has received little attention in previous studies.
Through this letter, we propose the analytically solvable non-Markovian non-Gaussian model. We also numerically reproduced the PDF of the real data which has heavy-tail. Thus we have shown that heavy-tail indeed is derived from memory.

{\it Acknowledgments}\,. - J.K. would like to thank to Gabjin Oh for his suggestions on the project, Jae-Sung Lee for the encouragement and advise, to Jung-Hyuck Park, Kanghoon Lee, Petre Jizba for the discussion and to Sang-Woo Lee and Sukjin Yun for the kind help in the simulation work. The work of T.O. is supported by Republic of Korea WCU grant R32-10130. J.K. and T.O. are thankful to Stephen Appleby for careful correction of the manuscript.  \\


\begin{thebibliography}{99}

\bibitem{Sornette}
D. Sornette, Critical Phenomena in Natural Sciences, Springer, Berlin, (2004).


\bibitem{deGennes}
P. de Gennes, Scaling concepts in polymer physics, SCornell University Press, (1979).


\bibitem{Peebles}
Peebles, P. J. E. ,The large-scale structure of the Universe, Princeton Series in Physics (1980).

\bibitem{PressSchechter}
Press, W. H. $\&$ Schechter, P., ApJ, 203, 557 (1976).

\bibitem{Kleinert}
 {\it Path Integrals in Quantum Mechanics, Statistics, Polymer Physics, and Financial Markets}\,, Hagen Keinert\,.
 

\bibitem{Cont}
R. CONT and P. TANKOV, Financial Modeling with Jump
Processes, Chapman and Hall/CRC (2004).

\bibitem{Rachev}
Rachev, S.T. and Mittnik, S. Stable Paretian Models in Finance. New York: Jone Wiley $\&$ Sons (2000).

 
\bibitem{Barabasi}
A.-L. Barabasi and R. Albert, Science 286, 509 (1999).



\bibitem{Ito}
Kiyoshi Ito, On stochastic differential equations. Memoirs, American Mathematical Society 4, 151 (1951).


\bibitem{Krug}
J. Franke, G. Wergen, and J. Krug, Correlations of Record Events as a Test for Heavy-Tailed Distributions, PRL 108, 064101 (2012).



\bibitem{Mandelbrot}
Mandelbrot, B.; van Ness, J.W. (1968), Fractional Brownian motions, fractional noises and applications, SIAM Review 10 (1968). 


\bibitem{Hurst}
 Hurst, H.E., Black, R.P., Simaika, Y.M., Long-term storage: an experimental study Constable, London (1965).


\bibitem{RSM}
Tong, H. and Lim, K. S.  Threshold Autoregression, Limit Cycles and Cyclical Data (with discussion), Journal of the Royal Statistical Society, Series B, 42, 245-292 (1980).

\bibitem{FriedmanUrn}
  A Simple Urn Model\,, Friedman, B., Comm. Pure Appl. Math.\, (1949) \,.

 
\end{thebibliography}
\end{document}